\documentclass[11pt,a4paper]{article}
\usepackage{amsmath}
\usepackage{amssymb}
\usepackage{amsthm}
\usepackage{latexsym}
\usepackage{amsbsy}
\usepackage[all]{xypic}
\usepackage{amscd}

\newtheorem{thm}{Theorem}[section]
\newtheorem{cor}[thm]{Corollary}
\newtheorem{lem}[thm]{Lemma}

\newtheorem{prop}[thm]{Proposition}
\newtheorem{defn}[thm]{Definition}
\newtheorem{rem}[thm]{Remark}

\begin{document}

\begin{center}
{\Large \bf Generalized Cluster Complexes via Quiver Representations}
\bigskip

{\large  Bin Zhu\footnote{Supported by the NSF of China (Grants
10471071) and by the Leverhulme Trust through the network 'Algebras,
Representations and Applications'.}}
\bigskip

{\small
 Department of Mathematical Sciences
\\  Tsinghua University
\\
    100084 Beijing, P. R. China
\\

{\footnotesize E-mail: bzhu@math.tsinghua.edu.cn}

}
\bigskip

\today

\end{center}

\def\s{\stackrel}
\def\gama{\gamma}
\def\Longrightarrow{{\longrightarrow}}
\def\P{{\cal P}}
\def\A{{\cal A}}
\def\F{\mathcal{F}}
\def\X{\mathcal{X}}
\def\T{\mathcal{T}}
\def\m{\textbf{ M}}
\def\t{{\tau }}
\def\b{\textbf{d}}
\def\K{{\cal K}}

\def\G{{\Gamma}}
\def\e{\mbox{exp}}

\def\righta{\rightarrow}

\def\s{\stackrel}

\def\ncong{\not\cong}

\def\mathbb{\NN}

\def\Hom{\mbox{Hom}}
\def\Ext{\mbox{Ext}}
\def\ind{\mbox{ind}}
\def\coprod{\amalg }
\def\L{\Lambda}
\def\c{\circ}
\def\mu{\multiput}

\renewcommand{\mod}{\operatorname{mod}\nolimits}
\newcommand{\add}{\operatorname{add}\nolimits}
\newcommand{\Rad}{\operatorname{Rad}\nolimits}
\newcommand{\RHom}{\operatorname{RHom}\nolimits}
\newcommand{\uHom}{\operatorname{\underline{Hom}}\nolimits}
\newcommand{\End}{\operatorname{End}\nolimits}
\renewcommand{\Im}{\operatorname{Im}\nolimits}
\newcommand{\Ker}{\operatorname{Ker}\nolimits}
\newcommand{\Coker}{\operatorname{Coker}\nolimits}
\renewcommand{\r}{\operatorname{\underline{r}}\nolimits}
\def \text{\mbox}

\hyphenation{ap-pro-xi-ma-tion}

\begin{abstract}
We give a quiver representation theoretic interpretation of
generalized cluster complexes defined by Fomin and Reading.
 By using $d-$cluster categories which are defined by Keller as triangulated orbit categories of (bounded)
  derived categories of
representations of valued quivers, we define a $d-$compatibility degree $(-||-)$ on any pair of ``colored''
almost positive real Schur roots
 which generalizes previous definitions on the non-colored case, and call two such roots
 compatible provided the
 $d-$compatibility degree of them is zero. Associated to the root system $\Phi$ corresponding to
 the valued quiver, by using this compatibility relation, we define
 a simplicial complex which has colored almost
positive real Schur roots as vertices and $d-$compatible subsets as
simplicies. If the valued quiver is an alternating quiver of a
Dynkin diagram, then this complex is the generalized cluster complex
 defined by Fomin and Reading.

\end{abstract}


\medskip


\medskip

\section{Introduction}
\medskip

Generalized cluster complexes associated to finite root systems are
introduced by Fomin and Reading \cite{FR2}. They have some nice
properties, see \cite{AT} and the references there. They are a
generalization of cluster complexes (so-called generalized
associahedra) associated to the same root systems introduced in
\cite{FZ2, FZ3}.
 Cluster complexes describe the combinatorial structure of cluster algebras which were
 introduced
by Fomin-Zelevinsky \cite{FZ1} in order to give an algebraic and
combinatorial framework for the canonical basis, see \cite{FR1} for
a nice survey on this combinatorics and also cluster combinatorics
of root systems. In \cite{MRZ}, Marsh-Reineke-Zelevinsky
 use ``decorated'' quiver
representations and tilting theory to give a quiver interpretation
of cluster complexes. This connection between tilting theory and
cluster combinatorics leads Buan-Marsh-Reineke-Reiten-Todorov
\cite{BMRRT} to introduce cluster categories for a categorical model
for cluster algebras, see also \cite{CCS} for type $A_n$. Cluster
categories are the orbit categories $\mathcal{D}/\tau^{-1}[1]$ of
derived categories of hereditary categories arising from the
 action of subgroup $<\tau^{-1}[1]>$ of the automorphism group. They are triangulated
categories \cite{Ke} and now they have become a successful model for
acyclic cluster algebras \cite{BMR, CC, CK}, see also the surveys
\cite{BM, Rin} and the references there for recent developments and
background of cluster tilting theory.

 $d-$cluster categories $\mathcal{D}/\tau^{-1}[d]$  as a generalization of cluster categories,
 were introduced by Keller \cite{Ke}, Thomas \cite{Th}, for $d\in \mathbf{N}$. They are studied by
 Keller and Reiten \cite{KR}, Y.Palu  \cite{Pa}, \cite{ABST}, see also \cite{BaM} for a geometric description of
 $d-$cluster categories of type $A_n$. $d-$cluster categories are triangulated categories
 with Calabi-Yau dimension $d+1$. When $d=1$, the cluster categories
 are recovered.

The aim of this paper is to give, not only a quiver representation
 theoretic interpretation of all key ingredients in
 defining generalized cluster complexes using $d-$cluster categories, but also a generalization
 of generalized cluster
 complexes to infinite root systems (compare Remark 3.13 in \cite{FR2}, there the authors asked whether there
 was such an extension).
  For simply-laced Dynkin case,
 Thomas \cite{Th} gives a realization of
generalized cluster complexes by defining the $d-$cluster categories.

The paper is organized as follows: In the first two parts, we recall
the well-known facts on $d-$cluster categories and (generalized)
cluster complexes of finite root systems. In particular, we recall
and generalize the BGP-reflection functors for cluster categories
\cite{Z1,Z2} to $d-$cluster categories. In the third part, we prove
 some properties of $d-$cluster tilting objects, including that any basic $d-$cluster tilting
 object contains exactly $n$ indecomposable direct summands. In the final section, for any
 root system $\Phi$, using a $d-$cluster category
$\mathcal{C}_d(\mathcal{H})$, we define a $d-$compatibility degree
on any pair of colored almost positive real Schur roots.  Using the
 $d-$compatibility degree, we define a generalized cluster complex
associated to $\Phi$, which has colored almost positive real Schur
roots as the vertices, and any subset forms a face if and
 only if any two elements of this subset are $d-$compatible. This simplicial complex is isomorphic to
the cluster complex of $d-$cluster category
$\mathcal{C}_d(\mathcal{H})$. If $\Phi$ is a finite root system, and
if we take $\mathcal{H}_0$ to be the category of representations of
an alternating quiver corresponding to $\Phi$, then our generalized
cluster complex is the usual generalized cluster complex
$\Delta^d(\Phi)$ defined
 by Fomin and Reading in \cite{FR2}.

\bigskip

\section{Basics on $d-$ cluster categories}

In this section, we collect some basic materials and fix notation which we will use later on.

 A valued graph $(\G,\b)$ is a finite set of
vertices $ 1,\cdots, n $,  together with non-negative integers
$d_{ij}$ for all pairs $i,j \in \G$ such that $d_{ii}=0$ and there
exist positive integers $\{\varepsilon _i\}_{i\in \G}$ satisfying
$$d_{ij}\varepsilon _j=d_{ji}\varepsilon _i,\  \ \mbox{for all }i,
\ j \in \G.$$  A pair $\{i, j\}$ of vertices is called an edge of
$(\G,\b)$ if
 $d_{ij}\not=0.$ An orientation $\Omega$ of a valued graph $(\G,\b)$
  is given by prescribing for each edge $\{i, j\}$ of $(\G,\b)$ an
  order (indicated by an arrow $i\rightarrow j$). For simplicity, we
  denote a valued graph by $\G$, and  a  valued quiver by $(\G,\Omega).$

Let $(\G,\Omega)$  be a valued quiver. We always assume that the
valued quiver $(\G, \Omega)$ contains no oriented cycles. Such
orientation $\Omega$ is called admissible.
  Let $K$ be a field and $\m=(F_i, {}_iM_j)_{i,j\in \G}$ a reduced $K-$species of
   $(\G,\Omega); $ that
  is, for all $i, j \in \G$, $_iM_j$ is an $F_i-F_j-$bimodule,
  where $F_i$ and $F_j$
   are division rings which are finite dimensional vector spaces over $K$
    and dim$(_{i}M_{j})_{F_j}=d_{ij}$
   and dim$_{K}F_i=
   \varepsilon_i$.
 We denote by $\mathcal{H}$ the category of finite dimensional representations of $(\G, \Omega, \mathcal{M})$.
 It is a hereditary abelian category \cite{DR}.
 Let $\Phi$ be the root system of the Kac-Moody Lie algebra
corresponding to the graph $\G$. We assume that
 $P_1, \cdots, P_n$ are the non-isomorphic indecomposable projective representations in $\mathcal{H}$, $E_1, \cdots, E_n$
 are the simple representations with dimension vectors $\alpha _1, \cdots, \alpha _n$,
 and  $\alpha_1, \cdots , \alpha_n$ are  the simple roots in $ \Phi.$ We
 use $D(-)$ to denote $\mbox{Hom}_K(-,K)$ which is a duality of $\mathcal{H}.$

Denote by $\mathcal{D} =
    D^{b}(\mathcal{H})$ the bounded derived category of $\mathcal{H}$ with shift
functor $[1]$.

\medskip

\subsection{$d-$cluster categories}

The derived category $\cal{D}$ has Auslander-Reiten triangles, and
the Auslander-Reiten translate $\tau$ is an automorphism of
$\mathcal{D}$. Fix a positive integer $d$, and denote by
$F_d=\tau^{-1}[d]$,
 it is an automorphism of $\mathcal{D}$. The $d-$cluster category of $H$ is defined in \cite{Ke, Th}:

We denote by $\mathcal{D}/ F_d$ the corresponding factor category.
The objects are by definition the $F_d$-orbits of objects in
$\cal{D}$, and the morphisms are given by
$$\Hom_{\mathcal{D}/F_d}(\widetilde{X},\widetilde{Y}) =
\oplus_{i \in \mathbf{Z}}
 \Hom_{\mathcal{D}}(X,F_d^iY).$$
Here $X$ and $Y$ are objects in $\cal{D}$, and $\widetilde{X}$ and
$\widetilde{Y}$ are the corresponding objects in $\mathcal{D}/F_d$
(although we shall sometimes write such objects simply as $X$ and
$Y$).

\begin{defn}\cite{Ke, Th}  The orbit category $\mathcal{D}/F_d$ is called the $d-$cluster category of $\mathcal{H}$
(or of $(\G,\Omega)$), which is denoted by
$\mathcal{C}_d(\mathcal{H})$, sometimes denoted  by
$\mathcal{C}_d(\Omega)$.\end{defn}

 By \cite{Ke}, the $d-$cluster category is a triangulated category with shift functor $[1]$ which
is induced by shift functor in $\mathcal{D}$, the projection
 $\pi: \mathcal{D}\longrightarrow \mathcal{D}/F$ is a triangle functor.
  When $d=1$, this orbit category is called the cluster category of $\mathcal{H}$, denoted by
  $\mathcal{C}(\mathcal{H})$ (sometimes denoted by $\mathcal{C}(\Omega ))$.

  $\mathcal{H}$ is a full subcategory of $\mathcal{D}$ consisting of complexes concentrated in
  degree $0$, then passing to $\mathcal{C}_d(\mathcal{H})$ by the projection $\pi$,  $\mathcal{H}$ is a
  (possibly not full) subcategory of $\mathcal{C}_d(\mathcal{H})$. For any $i\in \mathbf{Z}$, we use
  $(\mathcal{H})[i]$ to denote the copy of $\mathcal{H}$ under the $i-$th shift
$[i]$ as a subcategory of $\mathcal{C}_d(\mathcal{H}).$  In this
way, we have that $(\mbox{ind}\mathcal{H})[i]=\{ M[i]\ | \ M\in
\mbox{ind} \mathcal{H}\ \}$. For any object $M$ in
$\mathcal{C}_d(\mathcal{H})$, add$M$ denotes the full subcategory of
$\mathcal{C}_d(\mathcal{H})$ consisting of direct summands of direct
sums of copies of $M$.

For $X, Y\in \mathcal{C}_d(\mathcal{H})$, we will use Hom$(X,Y)$ to
denote the Hom-space Hom$_{\mathcal{C}_d(\mathcal{H})}(X,Y)$ in the
$d-$cluster category $\mathcal{C}_d(\mathcal{H})$ throughout the
paper. Define Ext$^i(X,Y)$ to be Hom$(X,Y[i]).$
\medskip

We summarize some known facts about $d-$cluster categories \cite{BMRRT, Ke}.

\begin{prop}\label{pr}\begin{enumerate}

\item $\mathcal{C}_d(\mathcal{H})$ has Auslander-Reiten triangles  and Serre functor
$\Sigma =\tau [1]$, where $\tau$ is the AR-translate in
$\mathcal{C}_d(\mathcal{H})$, which is induced from AR-translate in
$\mathcal{D}$.
\item   $\mathcal{C}_d(\mathcal{H})$ is a Calabi-Yau category of CY-dimension $d+1$.
\item  $\mathcal{C}_d(\mathcal{H})$ is  a Krull-Remark-Schmidt category.
\item $ \mathrm{ind }\mathcal{C}_d(\mathcal{H})=\bigcup_{i=0}^{i=d-1}(\mathrm{ind}\mathcal{H})[i]\bigcup\{ P_j[d]\ | \ 1\leq j\leq n\}.$
\end{enumerate}
\end{prop}

\begin{proof}\begin{enumerate}\item This is Proposition 1.3 \cite{BMRRT} and Corollary 1 in  Section 8.4 of \cite{Ke}.
\item It is proved in Corollary 1 in Section 8.4 of \cite{Ke}
\item This is proved in Proposition 1.2 \cite{BMRRT}
\item The proof for $d=1$ is given in Proposition 1.6 \cite{BMRRT}, which can be modified for the general $d$.
\end{enumerate}
\end{proof}

From Proposition \ref{pr}, we define the degree for every
indecomposable object in $\mathcal{C}_d(\mathcal{H})$ as follows:

\begin{defn}\label{degree} For any indecomposable object $X\in
\mathcal{C}_d(\mathcal{H})$, we call the non-negative integer
 $\mathrm{min }\{ k\in \mathbf{Z}_{\geq 0} \  | \  X \cong M[k] \mbox{ in } \mathcal{C}_d(\mathcal{H}),
 \mbox{ for  some  }  M\in \mathrm{ind }\mathcal{H}\ \}$ the degree of $X$, denoted by $\mathrm{deg }X$.
\end{defn}

 From Definition \ref{degree}, any indecomposable object $X$ of degree $k$ is isomorphic to $M[k]$ in
  $\mathcal{C}_d(\mathcal{H})$, where $M$ is an
 indecomposable representation in $\mathcal{H}$;  $0\leq \mbox{deg}X\leq d$, $X$ has degree $d$ if and
 only if $X\cong P[d]$ in $\mathcal{C}_d(\mathcal{H})$
 for some indecomposable projective object $P\in \mathcal{H}$; and $X$ has degree $0$ if and
 only if $X\cong M[0]$ in $\mathcal{C}_d(\mathcal{H})$
 for some indecomposable object $M\in \mathcal{H}$. Here $M[0]$ means regarding the object
 $M$ of $\mathcal{H}$
 as a complex concentrated in degree $0$.

\subsection{BGP-reflection functors}

If $T$ is a tilting object in $\mathcal{H}$, then the endomorphism
algebra $A=\mbox{End}_{\mathcal{H}}(T)$ is called a tilted algebra.
The tilting functor Hom$_{\mathcal{H}}(T,-)$ induces an equivalence
$\mbox{RHom}(T,-): D^b(\mathcal{H})\rightarrow D^b(A),$ where
$\mbox{RHom}(T,-)$ is the derived functor of
$\mbox{Hom}_{\mathcal{H}}(T,-).$

\medskip

Any standard triangle functor $G:  D^b(\mathcal{H})\rightarrow
D^b(\mathcal{H}')$ induces a well-defined functor
 $\tilde{G}: \mathcal{C}_d(\mathcal{H})\longrightarrow \mathcal{C}_d(\mathcal{H}')$ with
 the following commutative diagram \cite{Ke, Z1}:

\[ \begin{CD}
 D^b(\mathcal{H})  @>G>>  D^b(\mathcal{H}')\\
@V VV  @VV V  \\
\mathcal{C}_d(\mathcal{H})@>\tilde{G}
>> \mathcal{C}_d(\mathcal{H}')\end{CD} \]

The following result is proved in \cite{Z1, Z2}.
\begin{prop}
If $G:  D^b(\mathcal{H})\rightarrow D^b(\mathcal{H}')$ is a triangle
equivalence, then $\tilde{G }$ is also an equivalence of
triangulated categories.

\end{prop}

Let $k$ be a vertex in the valued quiver $(\G, \Omega )$, the
reflection of  $(\G, \Omega )$ at $k$ is the valued quiver  $(\G,
s_k\Omega )$, where $s_k\Omega$ is the orientation of $\G$ obtained
from $\Omega$ by reversing all arrows starting or ending at $k$. The
corresponding category of representations of $(\G,s_k\Omega,
\mathcal{M})$ is denoted simply by $s_k\mathcal{H}$. If $k$ is a
sink in the valued quiver $(\G, \Omega )$, then $k$ is a source of
$(\G, s_k\Omega )$, and the reflection of $(\G, s_k\Omega )$ at $k$
 is $(\G, \Omega )$. Let $k$ be a sink in $(\G, \Omega )$.
Then $P_k$ is a simple projective
 representation and $T=\oplus_{j\not=k} P_j\oplus \tau^{-1}P_k$ is a
tilting representation in $\mathcal{H}$ \cite{Rin}. The tilting
functor $S^+_k=\mbox{Hom}_{\mathcal{H}}(T,-)$ is a so-called
BGP-reflection functor, and its derived functor RHom$(T,-)$ is a
triangle
 equivalence from $D^b(\mathcal{H})$ to $D^b(s_k\mathcal{H})$, which is also denoted by $S^+_k.$  Similarly, one
has BGP-reflection functors $S^-_k$ for sources $k$.

\begin{defn} The induced functors
 $\widetilde{S^+_k} : \mathcal{C}_d(\mathcal{H})\longrightarrow \mathcal{C}_d(s_k\mathcal{H})$ for sinks $k$ and
$\widetilde{S^-_k} : \mathcal{C}_d(\mathcal{H})\longrightarrow
\mathcal{C}_d(s_k\mathcal{H})$ for sources $k$ are called
BGP-reflection functors of $d-$cluster categories.

\end{defn}

\begin{rem} When $d=1$, BGP-reflection functors are discussed in \cite{Z1}.\end{rem}
\medskip

We remind the reader that $\mathcal{H}$ (or $\mathcal{H}'$) is the
category of representations of the valued quiver $(\G, \Omega)$
($(\G, s_k\Omega)$, respectively);  the $P_i$ (respectively, the
$P'_i$)are the
 indecomposable projective representations in $\mathcal{H}$
(respectively, $\mathcal{H}'$) and the $E_i$ (respectively, the
$E'_i$) are the corresponding simple representations which are the
tops of the $P_i$ (respectively, the $P'_i$), for $i=1, \cdots , n.$

 \medskip

We recall from Proposition \ref{pr} and Definition \ref{degree} that
 any indecomposable object $Y$ in $\mathcal{C}_d(\mathcal{H})$ is isomorphic to $X[i]$ where $X\in
\mathrm{ind}\/\mathcal{H}$ and $i$ is the degree of $Y$. Keeping
this notation, we have the following proposition which gives the
images of indecomposable objects in
 $\mathcal{C}_d(\mathcal{H})$ under the BGP-reflection functor $\widetilde{S^+_k}.$

 \begin{prop} Let $k$
be a sink of the valued quiver $(\G, \Omega)$ and $Y$ an
indecomposable object in $\mathcal{C}_d(\mathcal{H})$ with degree
$i$.
 Then $Y\cong X[i]$ for an indecomposable representation $X$ in $\mathcal{H}$, and

$\widetilde{S^+_k}(X[i])= \left\{ \begin{array}{ll}
P_k'[d] &\mbox{ if  }X\cong P_k  (\cong E_k) \mbox{ and  } i=0 \\
 E_k'[i-1] &\mbox{ if  }X\cong P_k (\cong E_k) \mbox{ and  } 0<i\leq d \\
 P_j'[d] &\mbox{ if  }X\cong P_j \ncong P_k  \mbox{ and  } i=d \\
 S^+_k(X)[i] & \mbox{otherwise; }
\end{array}\right.$
     \end{prop}

\begin{proof} The statement in the proposition was proved in [Z1, Z2] when $d=1$. The proof for
the case $d>1$ is the same as there. We give a sketch of the proof
for the convenience of readers.  The BGP-reflection functor $S^+_k:
\mathcal{H}\longrightarrow s_k\mathcal{H}$ induces a triangle
equivalence $D^b(\mathcal{H})\longrightarrow D^b(s_k\mathcal{H}), $
denoted also by $S^+_k.$ It induces an equivalence
$\mathrm{ind}D^b(\mathcal{H})\longrightarrow
\mathrm{ind}D^b(s_k\mathcal{H})$. For any indecomposable object
$X[i]\in \mathrm{ind}D^b(\mathcal{H}),$ it is not hard to show that
$S^+_k(X[i])=S^+_k(X)[i]$ for $X\ncong P_k$ (note that $P_k=E_k $
since $k$ is a sink in $(\G,\Omega)$), and $S^+_k(P_k[i])=E'_k[i-1]$
for $i\in \textbf{Z}$ (compare [Z1] or [Z2]). Since $E'_k$ is an
injective representation in $s_k\mathcal{H}$, we have $\tau
P'_k[i]=E'_k[i-1]$ in $D^b(s_k\mathcal{H})$. Now passing to the
$d-$cluster category $\mathcal{C}_d(\mathcal{H})$ (which is an orbit
category of the derived category $D^b(\mathcal{H})$), we get the
images of indecomposable objects of $\mathcal{C}_d(\mathcal{H})$
under $\widetilde{S^+_k}$ as stated in the proposition.
\end{proof}
\medskip

\section{Cluster combinatorics of root systems}

For a valued graph $\G$, we denote by $\Phi=\Phi^+\bigcup\Phi^-$ the
set of roots of the corresponding Kac-Moody Lie algebra.

\begin{defn} \begin{enumerate} \item The set of almost positive roots is
 $$\Phi_{\geq -1}=\Phi^+\bigcup\{-\alpha_i\ | \ i=1,\cdots n\ \}.$$
\item Denote by $\Phi^{re}_{\geq-1}$ the subset of $\Phi_{\geq -1}$ consisting of
the positive real roots together with the negatives of the
simple roots.

\end{enumerate}

\end{defn}

When $\Phi$ is of finite type, $\Phi_{\geq -1}=\Phi^{re}_{\geq-1}$.

\begin{defn}   Let $s_i$ be the Coxeter generator of the
Weyl group of $\Phi$ corresponding to $i\in \G _0$. We call the following map
 the "truncated simple reflection" $ \sigma_i$ of $\Phi_{\geq
-1}$ \cite{FZ2}:
$$  \sigma_i(\alpha)=\left\{ \begin{array}{ll} \alpha & \alpha=-\alpha_j,
\ j\not=i \\ s_i(\alpha) & \mbox{otherwise.} \end{array}\right.$$
\end{defn}

It is easy to see that $\sigma_i$ is an automorphism of  $\Phi^{re}_{\geq-1}$.

\subsection{Cluster complexes of finite root systems}

In this first paragraph, we don't assume that $\G$ is a Dynkin
diagram (i.e. of finite type).
 Let $i_1,\cdots, i_n$ be an admissible ordering of $\G$ with respect to $\Omega$, i.e. $i_t$ is a sink with respect to
 $s_{i_{t-1}}\cdots s_{i_2}s_{i_1}\Omega$ for any $1\leq t\leq n.$
Denote by $R_{\Omega} =\sigma _{i_n}\cdots \sigma _{i_1}.$  This is
an automorphism of $\Phi _{\geq-1}$ and does not depend on the
choice of admissible ordering of $\G$ with respect to $\Omega.$ It
is the automorphism induced by Auslander-Reiten translation $\tau$
in $\mathcal{C}(\mathcal{H})$ (compare \cite{Z1, Z2}).

\medskip

In the rest of this subsection, we always assume that $\G$ is a
valued Dynkin graph, which is not necessarily connected.
 Fomin and Zelevinsky
[FZ2] associate a nonnegative integer $(\alpha||\beta)$, known as
the { compatibility degree}, to each pair $\alpha,\beta$ of almost
positive roots.

This is defined in the following way:  Let $\Omega _0$ denote one of
the alternating orientations of $\G$, and $\G ^+$ (respectively,
$\G^-$) the set of sinks (respectively, sources) of $(\G, \Omega_0)$
respectively. Define:
$$\ \ \ \ \tau_{\pm}=\prod_{i\in \G ^{\pm}}\sigma_i .$$

 Then
 $ R_{\Omega _0} =\tau_{-}\tau_{+}$, which is simply denoted
by $R$.
\medskip

 Denote by $n_i(\beta)$ the coefficient of $\alpha_i$ in the expansion of $\beta$ in terms of the simple roots
  $\alpha _1,
\cdots , \alpha _n$.
Then $(\ ||\ )$ is uniquely defined by the following two
properties:

 $$\begin{array}{lllcl}(*) & &(-\alpha_i||\beta)&=&\mbox{max}([\beta:\alpha _i],
 0),\\
  (**) &&(\tau_{\pm}\alpha||\tau_{\pm}\beta )&=&(\alpha ||\beta),\end{array}$$
  for any $\alpha , \beta \in \Phi_{\geq -1},$ any $i\in \G $.

\medskip

 Two almost positive roots $\alpha , \ \beta$  are called compatible if
$(\alpha ||\beta)=0$.

The cluster complex $\Delta (\Phi)$ associated to the finite root
system $\Phi$ is defined in \cite{FZ2}.

\begin{defn} The cluster complex $\Delta (\Phi)$ is a simplicial complex on the ground set $\Phi_{\geq
-1}$. Its faces are mutually compatible subsets of $\Phi_{\geq -1}.$
The facets of $\Delta (\Phi)$ are called the (root-)clusters
associated to $\Phi.$ \end{defn}

\subsection{Generalized cluster complexes of finite root systems}

 At the beginning of this subsection, we assume that $\G$ is
arbitrary valued graph, which is not necessarily connected, except
where we express specifically.
 As before $\Phi$ denotes the set of roots of the corresponding Lie algebra, and
 $\Phi_{\geq -1}$ denotes the set of almost positive roots. Fix a positive integer $d$, for any $\alpha \in \Phi ^+$,
 following \cite{FR2}, we call $\alpha ^1, \cdots, \alpha^d$ the $d$ ``colored'' copies of $\alpha$.

\begin{defn}\cite{FR2} The set of colored almost positive roots is
$$\Phi^d_{\geq-1}=\{ \alpha ^i: \alpha\in \Phi_{>0}, i\in \{1, \cdots, d\}\}\bigcup
\{(-\alpha _i)^1: 1\leq i\leq n\ \}.$$

\end{defn}

 When $\G$ is a Dynkin graph, the root system $\Phi$ of the corresponding Lie algebra is finite. In this case
 the generalized cluster complex $\Delta^d(\Phi)$ is defined on the ground set $\Phi^d_{\geq-1}$ and
using the binary compatibility relation on $\Phi^d_{\geq-1}$.
 This  binary compatibility relation is a natural generalization of  binary compatibility relation
 on $\Phi_{\geq-1}$ which we now recall from \cite{FR2}.

For a root $\beta\in \Phi_{\geq-1}$, let $t(\beta)$ denote the smallest $t$ such that $R ^t(\beta)$ is a
 negative root.

\begin{defn}\cite{FR2}  Two colored roots $\alpha ^k, \beta ^l\in \Phi^d_{\geq-1}$ are called compatible
if and only if one of the following conditions is satisfied:
\begin{enumerate}
\item $k>l.$ $t(\alpha)\leq t(\beta),$ and the roots $R(\alpha)$ and $\beta$ are compatible
(in the original ``non-colored'' sense)
\item  $k<l.$ $t(\alpha)\geq t(\beta),$ and the roots $\alpha$ and $R(\beta)$ are compatible
\item  $k>l.$ $t(\alpha)> t(\beta),$ and the roots $\alpha$ and $\beta$ are compatible
\item  $k<l.$ $t(\alpha) < t(\beta),$ and the roots $\alpha$ and $\beta$ are compatible
\item  $k=l.$ and the roots $\alpha$ and $\beta$ are compatible
\end{enumerate}

\end{defn}

Now we are ready to recall the definition of generalized cluster
complex  $\Delta^d (\Phi)$ for a finite root system $\Phi$.

\begin{defn} \cite{FR2} $\Delta^d (\Phi)$ has $\Phi^d_{\geq -1}$ as the set of vertices, its
simplices are mutually compatible subsets of $\Phi^d_{\geq -1}.$ The
subcomplex of  $\Delta^d (\Phi)$ which has $\Phi^d_{> 0}$ as the set
of vertices is denoted by  $\Delta^d_{+} (\Phi)$

\end{defn}

Now we generalize the definition of $R_d$ \cite{FR2} for finite root
system to an arbitrary root system.

\begin{defn} Let $( \G,\Omega)$ be a valued quiver. For $\alpha ^k\in \Phi^d_{\ge-1}, $ we set

$$R_{d,\Omega}(\alpha ^k)=\left\{
\begin{array}{ll}\alpha ^{k+1} & \mbox{ if }  \alpha \in \Phi_{>0} \mbox{ and } k<d;\\
&  \\
(R_{\Omega}(\alpha))^1 & \mbox{  otherwise  },\end{array}\right.$$

\end{defn}

\begin{rem} If $( \G,\Omega _0)$ is a valued Dynkin graph with an alternating orientation, then the automorphism
$R$ of $\Phi_{\geq -1}$ defined by Fomin-Zelevinsky \cite{FZ2} is
$R_{\Omega_0},$ hence $R_{d,\Omega _0}$ is the usual one ($R_d$)
defined by Fomin-Reading \cite{FR2}.\end{rem}

\begin{thm}\cite{FR2}  Let $\Phi$ be a finite root system. The compatibility relation on
 $\Phi^d_{\geq -1}$ has the following properties:
\begin{enumerate}\item $\alpha^k$ is compatible with $\beta^l$ if and only if $R_d(\alpha ^k)$ is
compatible with $R_d(\beta^l).$
\item $(-\alpha_i)^1$ is compatible with $\beta^l$ if and only if $n_i(\beta)=0.$

\end{enumerate}

Furthermore, conditions $1-2$ uniquely determine this relation.
\end{thm}

Now we generalize the ``truncated simple reflections'' of $\Phi_{\geq-1}$ to the colored almost positive roots.
Let $\Phi$ be an arbitrary root system (not necessarily of finite type).

\begin{defn}  Let $s_k$ be the Coxeter generator of the
Weyl group of $\Phi$ corresponding to $k\in \G_0 $. We define the following map
  $ \sigma_{k, d}$ of $\Phi^d_{\geq
-1}$:
$$  \sigma_{k, d}(\alpha^i)=\left\{ \begin{array}{ll} \alpha_k ^{d} & \mbox{ if  } i=1, \mbox{ and }\alpha=-\alpha_k,  \\
\alpha_k ^{i-1} & \mbox{ if  } 1<i\leq d, \mbox{ and } \alpha=\alpha_k, \\
 (-\alpha _j) ^1 & \mbox{ if  } i=1, \mbox{ and } \alpha=-\alpha_j,\ j\not=k
\\
(s_k(\alpha))^i & \mbox{otherwise.} \end{array}\right.$$
\end{defn}

$\sigma _{k,d}$ is a bijection of $\Phi ^d_{\geq -1}.$ We call it a
 $d-$truncated simple reflection of $\Phi^d_{\geq -1}.$

\medskip

\section{$d-$cluster tilting in $d-$cluster categories}

Let $\mathcal{C}_d(\mathcal{H})$ be a $d-$cluster category of type
$\mathcal{H}$, where $\mathcal{H}$ is the category of
representations of the valued quiver $(\G,\Omega)$. It is a
Calabi-Yau triangulated category with CY-dimension $ d+1 $.

\begin{defn}\label{exc}\begin{enumerate} \item  An object $X$ in
$\mathcal{C}_d(\mathcal{H})$ is called exceptional if
$\mathrm{Ext}^i(X,X)=0$, for any $1\leq i\leq d.$

\item An object $X$ is called a $d-$cluster tilting object if it satisfies the property: $Y\in
\mathrm{add}X$ if and only if $\mathrm{Ext}^i(X,Y)=0$ for $1\le i\le
d$.
\item An object $X$ is called almost complete tilting  if there is an indecomposable object
 $Y$ such that $X\oplus Y$ is a $d-$cluster tilting object.
  Such an indecomposable object $Y$ is called a complement of $X$.
\end{enumerate}
\end{defn}

\begin{prop}\label{basic} \begin{enumerate}\item For an object
$X$ in $\mathcal{H}$, $X$ is exceptional in $\mathcal{H}$ i.e.
 $\mathrm{Ext}^1_{\mathcal{H}}(X,X)=0$  if and only if $X[0]$ is
exceptional in $\mathcal{C}_d (\mathcal{H}).$
\item Any indecomposable exceptional object $X$ in $\mathcal{C}_d(\mathcal{H})$
is of the form $M[i]$ with $M$ being an exceptional representation
in $\mathcal{H}$ and $0\leq i\leq d-1$ or of the form $P_j[d]$ for
some $1\leq j\leq n$. In particular, if $\G$ is a Dynkin graph, then
any indecomposable object in $\mathcal{C}_d(\mathcal{H})$ is
exceptional.

\item Suppose $d>1$. Then $\mathrm{End}_{\mathcal{C}_d(\mathcal{H})}X$ is a division algebra
 for any indecomposable exceptional object $X.$
\item Suppose $d>1$. Let $P$ be a projective representation in $\mathcal{H}$ and $X$
a representation in $\mathcal{H}$. Then for any $-d\leq i\le d$,
$\mathrm{Ext}^1(P,X[i])=0$ except possibly for $i\in \{-1,\ d-1,\ d
\}$.
\end{enumerate}

\end{prop}
\begin{proof}\begin{enumerate}\item Let $X\in \mathcal{H}$ be exceptional. We will prove that $\mbox{Ext}^i(X,X)=0,$
for any $i\in \{ 1, \cdots, d\}.$ By definition, we have that
$\mbox{Ext}^i(X,X)=\oplus _{k\in
\mathbf{Z}}\mbox{Ext}^i_{\mathcal{D}}(X,
\tau^{-k}X[kd])=\mbox{Ext}^i_{\mathcal{D}}(X,X)\oplus
\mbox{Ext}^i_{\mathcal{D}}(X,\tau X[-d])$. In this sum, the first
summand $\mbox{Ext}^i_{\mathcal{D}}(X,X)=0, \forall i\ge 1,$ while
the second summand $\mbox{Ext}^i_{\mathcal{D}}(X,\tau X[-d])\cong
\mbox{Hom}_{\mathcal{D}}(X, \tau X[i-d]),$ which is zero when $i<d,$
and is isomorphic to $\mbox{Ext}^1_{\mathcal{D}}(X, X)=0$  when
$i=d.$ This proves that $X$ is exceptional in
$\mathcal{C}_d(\mathcal{H}).$ The proof for the converse follows
directly from definition.

\item The statements follow from Proposition 2.2(4) and Definition \ref{exc}, also using part 1 and the fact that
the shift is an autoequivalence.

\item Let $X$ be an indecomposable exceptional representation in $\mathcal{H}$ and suppose $d>1.$  It follows from
the definition of the orbit category that
End$_{\mathcal{C}_d(H)}X\cong \oplus_{m\in
\mathbf{Z}}\mbox{Hom}_{\mathcal{D}}(X,\tau^{-m}X[dm])\\
  \cong
\mbox{End}_{\mathcal{H}}X$. The last isomorphism holds due to the
facts: $\mbox{Hom}_{\mathcal{D}}(X, \tau ^{m}X[-md])\cong
\mbox{Hom}_{\mathcal{D}}(X[md], \tau ^{m}X)\cong
\mbox{Ext}^1_{\mathcal{D}}(\tau^{m-1}X, X[md]) =0$ for any positive
integer $m$; and $\mbox{Hom}_{\mathcal{D}}(X, \tau ^{-m}X[md])\cong
\mbox{Hom}_{\mathcal{D}}(\tau ^mX, X[md])$, which is also zero since
 $md>1$ (we use the assumption $d>1$ here) for any positive integer $m$.
 Then End$_{\mathcal{C}_d(H)}X $ is a division algebra since
 End$_{\mathcal{H}}X$ is a division algebra. Since any indecomposable exceptional object $M$ in $\mathcal{C}_d(H)$ is
some shift $X[i]$ of an indecomposable exceptional representation
$X$ in $\mathcal{H}, $
 $\mbox{End}_{\mathcal{C}_d(H)}M= \mbox{End}_{\mathcal{C}_d(H)}X[i]\cong \mbox{End}_{\mathcal{C}_d(H)}X$
 is a division algebra.
\item Suppose $d>1$. Let $P$ be a projective representation in $\mathcal{H}$ and $X$ a representation
in $\mathcal{H}$. Then for any $-d\le i\le d$,
$\mbox{Ext}^1(P,X[i])=\oplus _{k\in
\mathbf{Z}}\mbox{Ext}^1_{\mathcal{D}}(P, \tau^{-k}X[dk+i]) \cong
\mbox{Ext}^1_{\mathcal{D}}(P, \tau X[-d+i])\oplus
\mbox{Ext}^1_{\mathcal{D}}(P,X[i])$. Now if $i\not=-1,\ d-1 , \ d,$
then $\mbox{Ext}^1_{\mathcal{D}}(P, \tau X[-d+i])=0=
\mbox{Ext}^1_{\mathcal{D}}(P,X[i])$.
 Then for any $-d\le i\le d$, $\mbox{Ext}^1(P,X[i])=0$ except for $i=-1,\   d-1$ and $d$.
\end{enumerate}
\end{proof}

\begin{rem}\label{rem} Any basic (i.e. multiplicity-free) exceptional object contains at most $(d+1)n$ non-isomorphic indecomposable direct summands. \end{rem}

\begin{proof} Let $X$ be a basic exceptional object in $\mathcal{C}_d(\mathcal{H})$. Then any indecomposable direct summand of $X$ is exceptional, hence by  Proposition \ref{basic} (2), we write $M$ as $M=\oplus _{k=0}^{k=d}\oplus_{i\in I_k}M_{i,k}[k]$ with $M_{i, k}$ being an indecomposable exceptional representation. Therefore $\oplus _{i\in I_k}M_{i,k}$ is an exceptional object in hereditary category $\mathcal{H}$, hence the number of direct summands is at most $n$, i.e. $|I_k|\leq n$. Then the number of indecomposable direct summands of $M$ is at most $(d+1)n$.
\end{proof}

For any pair of objects $T,\  X$ in $\mathcal{C}_d(\mathcal{H})$,
due to the Calabi-Yau property of $\mathcal{C}_d(\mathcal{H})$, we
have that $\mathrm{Ext}^i(X,T)=0$ for $1\le i\le d$ if and only if
$\mathrm{Ext}^i(T,X)=0$ for $1\le i\le d$. Hence from Remark
\ref{rem} and Definition \ref{exc}, $T$ is a $d-$cluster tilting
object in $\mathcal{C}_d(\mathcal{H})$ if and only if add$T$ is a
maximal $d-$orthogonal subcategory of $\mathcal{C}_d(\mathcal{H})$
in the sense of \cite{I2}: i.e. add$T$ is contravariantly finite and
covariantly finite in $\mathcal{C}_d(\mathcal{H})$, and satisfies
the property: $X\in \mathrm{add}T$ if and only if
$\mathrm{Ext}^i(X,T)=0$ for $1\le i\le d$ if and only if
$\mathrm{Ext}^i(T,X)=0$ for $1\le i\le d$. In the following we will
prove that any basic $d-$cluster tilting object contains exactly $n$
indecomposable direct summands. First of all, we recall some results
from \cite{I2} which hold in any $(d+1)-$Calabi-Yau triangulated
category.

\begin{thm} (Iyama)\label{Iyama}  Let $X$ be an almost complete tilting object in $\mathcal{C}_d(\mathcal{H})$, and $X_0$ a complement of $X$. Then there are $d+1$ triangles:
$$(*)\  \  \  X_{i+1}\s{g_i}{\longrightarrow}B_i\s{f_i}{\longrightarrow}X_i\s{\sigma _i}{\longrightarrow}X_{i+1}[1]$$
where $f_i$ is the minimal right $\mathrm{add}{X}-$approximation of
$X_i$ and $g_i$ minimal left $\mathrm{add}{X}-$approximation of
$X_{i+1}$, all $X_i$ are indecomposable and complements of $X$,
$i=0, \cdots, d$.
\end{thm}

For the convenience of readers, we sketch the proof, for details,
see \cite{IY}.

\begin{proof} We suppose that $d>1$,  the same statement for $d=1$ was proved in \cite{BMRRT}.
 For the complement $X_0$ of $X$, we
consider the minimal right add$X-$approximation $f_0: B_0\rightarrow
X_0$ of $X_0$, extend $f_0$ to a triangle  $X_1\s{g_0}{\rightarrow}
B_0\s{f_0}{\rightarrow}X_0\s{\sigma _0}{\rightarrow} X_1[1].$ It is
easy to see that $X_1$ is indecomposable, $g_0$ is the minimal left
add$X-$approximation of $X_1$ and $X\oplus X_1$ is an exceptional
object in $\mathcal{C}_d(\mathcal{H})$ (compare \cite{BMRRT}). It
follows from Theorem 5.1 in \cite{IY} that $X\oplus X_1$ is a
$d-$cluster tilting object. Continuing this step, one can get
complements $X_1, \cdots, X_{d+1}$ with triangles
$X_{i+1}\s{g_{i}}{\rightarrow}
B_{i}\s{f_{i}}{\rightarrow}X_{i}\s{\sigma _i}{\rightarrow} X_{i}[1]$
for $0\le i\le d$, where $f_i$ ($g_i$) is the minimal right (left,
resp.) add$X-$approximation of $X_i$ $(X_{i+1},$ resp.) and $X\oplus
X_i$ is a $d-$cluster tilting object.

\end{proof}

\begin{cor} \label{Iyama2}  With the same notation as Theorem \ref{Iyama},
we have that $\sigma_d[d]\sigma_{d-1}[d-1]\cdots
\sigma_1[1]\sigma_0\not=0.$ In particular,
$\mathrm{Hom}(X_i,X_j[j-i])\not=0$ and $X_i\ncong X_j, \forall 0\leq
i<j\leq d.$
\end{cor}
\begin{proof}
From Theorem \ref{Iyama}, we have that $\sigma_0\not=0$ since the
 triangle $(*)$ at $i=0$ in Theorem \ref{Iyama} is non-splitting. Suppose that
$\sigma_d[d]\sigma_{d-1}[d-1]\cdots \sigma_1[1]\sigma_0=0$, then
$\sigma_{d-1}[d-1]\cdots \sigma_1[1]\sigma_0: X_0\rightarrow
X_{d}[d]$ factors through $f_d[d]: B_{d}[d]\rightarrow X_{d}[d]$
since we have a triangle
$X_{d+1}[d]\s{g_d[d]}{\longrightarrow}B_{d}[d]\s{f_d[d]}{\longrightarrow}X_{d}[d]\s{\sigma
_d[d]}{\longrightarrow}X_{d+1}[d+1]$. Since Hom$(X_0,
B_d[d])=\mbox{Ext}^d(X_0,B_d)=0$, $\sigma_{d-1}[d-1]\cdots
\sigma_1[1]\sigma_0=0.$
 Similarly $\sigma_{d-2}[d-2]\cdots \sigma_1[1]\sigma_0=0$, and finally, $\sigma_0=0$, a contradiction.
 Now we prove final statement: we have that
 $\sigma_{j-1}[j-1]\cdots \sigma_i\in \mbox{Hom}(X_i, X_j[j-i])$ and
 $\sigma_{j-1}[j-1]\cdots \sigma_i\not= 0.$ Otherwise $\sigma_{j-1}[j-1]\cdots \sigma_i= 0$, and
 hence $\sigma_{d-1}[d-1]\cdots \sigma_1[1]\sigma_0=0,$ a
 contradiction.  Now suppose that $X_i\cong X_j$ for some $i<j$. Then Ext$^k(X_i,X_j)=0$ for $1\leq k \leq d$,
 a contradiction. Then $X_i\ncong X_j$.
\end{proof}

Now we state our main result of this section.

\begin{thm}\label{tilting}

 Any basic $d-$cluster tilting object in $\mathcal{C}_d(\mathcal{H})$ contains exactly $n$ indecomposable direct summands.

\end{thm}

To prove the theorem, we need some technical lemmas.

\begin{lem}\label{deg} Let $d>1$ and  $X=M[i], \ Y=N[j]$ be indecomposable objects of degree $i, \ j$
respectively in $\mathcal{C}_d(\mathcal{H})$. Suppose that
$\mathrm{Hom}(X, Y)\not=0.$ Then one of the following holds: (1) We
have $i=j$ or $j-1$ (provided $j\geq1$);

(2) We have $i=0$, $i=d$  (and  $M=P)$ or $d-1$ (provided $j=0)$.
\end{lem}

\begin{proof} Let $d>1$. Firstly we note that for any indecomposable object $X\in \mathcal{C}_d(\mathcal{H})$,
  $0\leq \mbox{deg}X\leq d$, deg$X=d$ if and only if $X=P_i[d]$ for an indecomposable projective representation $P_i$.
      This implies that $-d\leq \mbox{deg}Y-\mbox{deg}X \leq d$ for indecomposable objects $X, Y\in \mathcal{C}_d(\mathcal{H})$.
       Let $X=M[i], \ Y=N[j]$ be indecomposable objects of degree $i, \ j$ respectively in
$\mathcal{C}_d(\mathcal{H})$. We have $\mbox{Hom}(X,Y)\cong
\mbox{Hom}(M,N[j-i])=\oplus_{k\in
\mathbf{Z}}\mbox{Hom}_{\mathcal{D}}(M, \tau^{-k}N[j-i+kd])
 =\mbox{Hom}_{\mathcal{D}}(M, \tau N[j-i-d])\oplus \mbox{Hom}_{\mathcal{D}}(M,
N[j-i]) \oplus \mbox{Hom}_{\mathcal{D}}(M, \tau^{-1} N[j-i+d])$. The
last equality holds  due to $-d\le j-i\le d$, and
Hom$_{\mathcal{D}}(M, \tau^{-k}N[j-i+kd])=0$ for $k\not=-1, 0, 1$.
We divide the calculation of $\mbox{Hom}(X,Y)$ into three cases:
\begin{enumerate}\item The case $-d<j-i<d.$  We have that $\mbox{Hom}(X, Y)\cong
\mbox{Hom}_{\mathcal{D}}(M, N[j-i])\oplus
\mbox{Hom}_{\mathcal{D}}(M, \tau ^{-1}N[j-i+d])$. The first summand
is zero when $j-i\not= 0, 1$
 while the second is zero when  $d+j-i\not= 1$ (equivalently $d+j-i>1$ since $0<d+j-i<2d$).
\item The case $j-i=-d$. Then $j=0, \ i=d \ (M=P).$ Then $\mbox{Hom}(X, Y)=\mbox{Hom}_{\mathcal{D}}(P, \tau ^{-1}N)$.
 \item The case $j-i=d.$ Then $j=d \ (N=P), i=0.$ Then $\mbox{Hom}(X, Y)
 =\mbox{Hom}_{\mathcal{D}}(M, \tau P)\oplus  \mbox{Hom}_{\mathcal{D}}(M, P[d])\oplus
 \mbox{Hom}_{\mathcal{D}}(M, \tau^{-1} P[2d])=0$.
\end{enumerate}
Therefore if $\mbox{Hom}(X, Y)\not=0$, then
$\mbox{Hom}(M[0],N[j-i])\not=0$. Proof of (1). Suppose $j\geq 1$.
Then combining with Case $3,$ we have that $-d<j-i<d.$ We want to
prove that if $j-i\not= 0,1$ then $\mbox{Hom}(X, Y)=0$ this will
finish the proof of $(1)$. Under the condition $j-i\not= 0,1$, from
Case 1, we have that $\mbox{Hom}(X, Y)\cong
\mbox{Hom}_{\mathcal{D}}(M, \tau ^{-1}N[j-i+d])$, which is zero for
$d+j-i\not= 1$. But if $d+j-i= 1$, i.e. $i=d,$ hence $M=P$ and
$j=1$. Then $\mbox{Hom}_{\mathcal{D}}(M, \tau
^{-1}N[j-i+d])=\mbox{Hom}_{\mathcal{D}}(P, \tau ^{-1}N[1])=0$. We
 have finished the proof of (1). Proof of (2). Suppose $j=0.$  Then $-d\leq j-i\leq 0$.
It follows from Cases 1-2 that $i=0, \ i=d \ (M=P)$ or $i=d-1.$ This
finishes the proof of (2).

\end{proof}

\begin{lem}\label{ext} If $d>2$, then $\mathrm{Ext}^2(M[i],N[i])=0$ for objects $M, N\in \mathcal{H}$, and any $i$.
\end{lem}
\begin{proof}
It is sufficient to prove that $\mbox{Ext}^2(M[0],N[0])=0$. From the
definition of the orbit category $\mathcal{D}/\tau^{-1}[d]$, we have
that
 $$\mbox{Ext}^2(M[0],N[0])=\mbox{Hom}(M[0],N[2])=\oplus _{k\in \mathbf{Z}}\mbox{Hom}_{\mathcal{D}}(M, \tau^{-k}N[kd+2]),$$
 where each summand
 $\mbox{Hom}_{\mathcal{D}}(M, \tau^{-k}N[kd+2])$ equals $0$ since $kd+2\ge 2$ or $kd+2\leq -1$ by the condition $d>2.$
 Hence $\mbox{Ext}^2(M[0],N[0])=0.$
\end{proof}

\begin{lem}\label{d=2} Let $d>1$ and $M, N\in \mathcal{H}$. Then $\mathrm{Ext}^1(M[0],N[0])\cong
\mathrm{Ext}^1_{\mathcal{H}}(M,N)$. Furthermore any non-split
triangle between $M[0]$ and $N[0]$ in $\mathcal{C}_d(\mathcal{H})$
is induced from a non-split exact sequence between $M$ and $N$ in
$\mathcal{H}$.
\end{lem}

\begin{proof} Under the condition $d>1$, it is easy to see that
$\mbox{Ext}^1(M[0],N[0])=\\ \oplus _{k\in \mathbf{Z}}
\mbox{Ext}^1_{\mathcal{D}}(M, \tau^{-k}N[2k])=
\mbox{Ext}^1_{\mathcal{D}}(M,N) =\mbox{Ext}^1_{\mathcal{H}}(M,N).$
This proves the first statement. Since $\mathcal{H}\subset
\mathcal{C}_d(\mathcal{H})$ is a (not necessarily full) embedding,
 and any exact short sequence in $\mathcal{H}$ induces a triangle in
$ \mathcal{C}_d(\mathcal{H})$, the final statement then follows from
the first statement.

\end{proof}

\begin{proof}(of Theorem \ref{tilting}):
 We assume that $d>1$ since it was proved in [BMMRT] for $d=1$.
  Let $M=\oplus_{i\in I} M_i[k_i]$ be a $d-$cluster tilting object in
   $\mathcal{C}_d(\mathcal{H})$, where all $M_i$ are indecomposable representations in $\mathcal{H}$,
$0\le k_i\le d$ (when $k_i=d$, $M_i$ is projective). One can assume
that one of $k_i$ is $0$, otherwise one can replace $M$ by a
suitable shift of $M$.  Denote by $\nu (M)=max\{|k_i-k_j|\  |\
\forall i, j\}.$  We prove $|I|=n$ by induction on $\nu(M)$, where
$|I|$ denotes the cardinality of $I$.  If $\nu(M)=0$, i.e. $k_i=0$
for all $i$, then $\oplus_{i\in I} M_i[0]$ is a $d-$cluster tilting
object in $\mathcal{C}_d(\mathcal{H})$, hence a tilting object in
$\mathcal{H}.$  Then $|I|=n.$ Now assume  $\nu(M)=m>0.$ Without
losing generality, we assume that $k_1=\cdots =k_t=m$ and $k_j<m$
for $j>t$.  From the complement $X_0=M_1[k_1]$ of $X=M\setminus
M_1[k_1]$ (here we use $X \setminus X_1$ to denote a complement of
$X_1$ in $X$ for a direct summand $X_1$ of $X$), by Theorem
\ref{Iyama}, we have at least $d+1$ complements $X_j$, $j=0, \cdots,
d,$ which form the triangles $(*)$ in  Theorem \ref{Iyama}. In these
triangles, it is easy to see that $f_i= 0$ if and only if $B_i=0$ if
and only if $g_i= 0$. We will prove that there are at least one of
complements $X_j$ with smaller degree than $m$. At first, we prove
this statement for the special case that $m=1$. We claim that the
degree of $X_1$ is $0$ or $1$ in this case. Otherwise $X_1=P[d]$ for
some indecomposable projective representation $P$ or $X_1=Y[d-1]$
for some indecomposable representation $Y$. Write $X_0$ as $Z[1]$,
where $Z$ is an indecomposable representation in $\mathcal{H}$.
 If $X_1=P[d]$, then
Hom$(X_1,X_0[d])=\mbox{Hom}(P[d],X_0[d])\cong\mbox{Hom}(P,Z[1])=0$,
this contradicts to the fact that Hom$(X_1,X_0[d])\cong
\mbox{Hom}(X_0, X_1[1])$ is not zero by Theorem 4.4 or Corollary
4.5.  If $X_1=Y[d-1]$, then $X_1$ has degree $1$ when $d=2$, and
Hom$(X_1,X_0[d])=\mbox{Hom}(Y[d-1],Z[d+1])\cong\mbox{Ext}^2(Y,Z)=0$
by Lemma 4.8. when $d>2$, which also contradicts  to the fact that
Hom$(X_1,X_0[d])\cong \mbox{Hom}(X_0, X_1[1])$ is not zero. This
proves the statement that $X_1$ has degree $0$ or $1$.  Now if there
are no complements $X_j$ of $X$ with degree $0$, then all $X_j$ have
degree $1$. We prove that any three successive complements, say
$X_0, X_1, X_2$, can not have the same degree. If all degrees of
$X_i,i=0,1,2,$ are the same, we can assume that all $X_i$ have
degree $0$. By Lemma \ref{d=2}, we have non-split short exact
sequences in $\mathcal{H}$:
$$\begin{array}{l}0\longrightarrow X_1\longrightarrow B_0 \longrightarrow X_0\longrightarrow 0\\
0\longrightarrow X_2\longrightarrow B_1 \longrightarrow
X_1\longrightarrow 0\end{array}$$
 From the first short exact sequence, we have $\mbox{Ext}^1_{\mathcal{H}}(X_0,X_1)\ncong 0$. By applying
 $\mbox{Hom}_{\mathcal{H}}(X_0,-)$ to the second exact sequence, we have the
 long exact sequence:
 $$\cdots\rightarrow \mathrm{Ext}^1_{\mathcal{H}}(X_0,X_2)\rightarrow \mathrm{Ext}^1_{\mathcal{H}}(X_0,B_1) \rightarrow
 \mathrm{Ext}^1_{\mathcal{H}}(X_0,X_1)\rightarrow \mathrm{Ext}^2_{\mathcal{H}}(X_0,X_2)\rightarrow
 \mathrm{Ext}^2_{\mathcal{H}}(X_0,B_1).$$  Since $X\oplus X_0$ is a $d-$cluster tilting object
 in $\mathcal{C}_d(\mathcal{H})$ and $B_1\in \mbox{add} X$,  $\mathrm{Ext}^1(X_0,B_1)=0$.
 Hence we have that
 $\mathrm{Ext}^1_{\mathcal{H}}(X_0,B_1)=0$ by Lemma 4.9.
 It follows that $\mbox{Ext}^1_{\mathcal{H}}(X_0,X_1)= 0$ since
 $\mbox{Ext}^2_{\mathcal{H}}(X_0,X_2)=0$ due to $\mathcal{H}$ being hereditary. It is a
 contradiction. This finishes the proof for $m=1$.

 Now suppose $m>1$. We will prove that there are at least one of complements $X_j$ with
smaller degree than $m$.
 We divide the proof into two cases: Case 1. All maps $f_i$
(equivalently $g_i$) are non-zero.  Now we assume that there are no
complements of $X$ with smaller degree than $m$. Then by Lemma 4.7,
the degrees of all $X_i$ are $m$. If $d>2$,
 then $\mbox{Ext}^2(X_0, X_2)=0$ by Lemma \ref{ext}, which is a
contradiction to Corollary \ref{Iyama2}.
 If $d=2$, then the same proof as above shows
  that $\mbox{Ext}^1(X_0,X_1)= 0$ which contradicts to Corollary \ref{Iyama2}.
  Therefore there is a complement of $X$ with smaller degree than $m$.
 Case 2. There are some $i$ such that $f_i=0$ (equivalently $g_i=0$).  Then $X_{i}\cong
X_{i+1}[1]$ for such $i$. It follows that $X_{i+1}$ has smaller
degree than $X_i$ if $X_i$ has strictly positive degree.

 Therefore we have a complement of $X$, say $X_{s}$, such that the degree ${k'_1}$ of $X_s$ is smaller than $m=k_1$.
 Now we replace $X$ by $X'=(X\setminus X_0)
 \oplus X_s$, which is, by Theorem \ref{Iyama},  a $d-$cluster  tilting object in $\mathcal{C}_d(\mathcal{H})$, containing $|I|$ indecomposable direct summands. The number of indecomposable
direct summands of $X'$ with the (maximal) degree $m(=\nu (M))$ is
$t-1$. We repeat the step for the complement $M_2[k_2]$ of almost
complete tilting object $X'\setminus M_2[k_2]$, we get a $d-$cluster
tilting object $X''$ containing $|I|$ indecomposable direct
summands, and the number of indecomposable direct summands of $X''$
with the (maximal) degree $m(=\nu (M))$ is $t-2$.  Repeating such
step $t$ times, one can get a (basic) $d-$cluster tilting object $T$
containing $|I|$ indecomposable direct summands and $\nu(T)<\nu
(M).$  By induction, $T$ contains exactly $n$ indecomposable direct
summands. Then $|I|=n.$

\end{proof}

\begin{rem} Theorem \ref{tilting} is proved by Thomas in \cite{Th} for a simply-laced Dynkin
quiver $(\G,\Omega _0)$, there he uses the fact that
 $\mathrm{ind}D^b(K\vec{\Delta})\approx \mathbf{Z}\vec{\Delta}$ for a Dynkin quiver $\vec{\Delta}$. This fact does not hold for
 non-Dynkin quivers. Our proof is more categorical.\end{rem}

Denote by $\mathcal{E}(\mathcal{H})$ the set of isomorphism classes
of indecomposable exceptional  representations in $\mathcal{H}$. The
set $\mathcal{E}(\mathcal{C}_d(\mathcal{H}))$ of isoclasses of
indecomposable exceptional objects in $\mathcal{C}_d(\mathcal{H})$
is the (disjoint) union of subsets  $\mathcal{E}(\mathcal{H})[i]$,
$i=0,1, \cdots, d-1,$ with $\{ P_j[d]| 1\leq j\leq n\}$. A subset
$\mathcal{M}$ of $\mathcal{E}(\mathcal{C}_d(\mathcal{H}))$ is called
exceptional if for any $X, Y\in \mathcal{M}$,
$\mathrm{Ext}^i(X,Y)=0$ for all $i=1,\cdots, d$. Denote by
$\mathcal{E}_+(\mathcal{C}_d(\mathcal{H}))$ the subset of
$\mathcal{E}(\mathcal{C}_d(\mathcal{H}))$ consisting of all
indecomposable exceptional objects other than
 $P_1[d], \cdots P_n[d].$

Now we are ready to define a simplicial complex associated to the
$d-$cluster category $\mathcal{C}_d(\mathcal{H})$, which is a
generalization of the classical cluster complexes of cluster
categories \cite{BMRRT, Rin, Z1}.

\begin{defn} The cluster complex $\Delta ^d(\mathcal{H})$ of $\mathcal{C}_d(\mathcal{H})$
is a simplicial complex which has
$\mathcal{E}(\mathcal{C}_d(\mathcal{H}))$ as the set of vertices,
and has exceptional subsets in $\mathcal{C}_d(\mathcal{H})$ as its
simplices. The positive part $\Delta ^d_{+}(\mathcal{H})$ is the
subcomplex of $\Delta ^d(\mathcal{H})$ on the subset
$\mathcal{E}_+(\mathcal{C}_d(\mathcal{H}))$.

\end{defn}

From the definition,  the facets (maximal simplices) are exactly the
$d-$cluster tilting subsets (i.e. the sets of indecomposable objects
of $\mathcal{C}_d(\mathcal{H})$ (up to isomorphism) whose direct sum
 is a $d-$cluster tilting object).

\begin{prop}\begin{enumerate}\item  $\Delta ^d(\mathcal{H})$ and  $\Delta ^d_{+}(\mathcal{H})$ are pure of dimension $n-1$,
\item  For any sink (or source) $k$, the BGP-reflection functor $\tilde{S_k^+}$ (resp.  $\tilde{S_k^-}$) induces an isomorphism between
$\Delta ^d(\mathcal{H})$ and $\Delta ^d(s_k\mathcal{H})$. In
particular, if $\G$ is a Dynkin diagram and $\Omega$ and $\Omega'$
are two orientations of $\G$, then $\Delta ^d(\mathcal{H})$ and
$\Delta ^d(\mathcal{H}')$ are isomorphic.

\end{enumerate}
\end{prop}

\begin{proof}\begin{enumerate}\item It follows from Theorem \ref{tilting} that any
$d-$cluster tilting subset contains exactly $n$ elements. Hence
$\Delta ^d(\mathcal{H})$ is pure of dimension $n-1$. Now suppose
$M=\oplus_{i=1}^{n-1} M_i$ is an exceptional object
 in $\mathcal{C}_d(\mathcal{H})$
and none of the $M_i$ are isomorphic to $P_j[d]$ for any $j$. In the
proof of Theorem \ref{tilting}, we proved that not all complements
of an almost complete tilting objects have the same degrees. Then
$M$ has a complement in $\mathcal{E}_+(\mathcal{C}_d(\mathcal{H})$.
This proves that
 $\Delta ^d_{+}(\mathcal{H})$ is pure of dimension $n-1$.
\item Since $\tilde{S_k^+}$ is a triangle equivalence from the $d-$cluster
category $\mathcal{C}_d(\mathcal{H})$ to
$\mathcal{C}_d(s_k\mathcal{H})$, it sends (indecomposable)
 exceptional objects to (indecomposable)
exceptional objects. Thus it induces an isomorphism from $\Delta
^d(\mathcal{H})$ to $\Delta ^d(s_k\mathcal{H})$. The second
statement follows from the first statement together with
 that for two orientations $\Omega,
\ \Omega '$ of a Dynkin graph $\G$, there is a admissible sequence
with respect to sinks $i_1, \cdots, i_n$ such that $\Omega
'=s_{i_n}\cdots s_{i_1}\Omega.$
\end{enumerate}
\end{proof}
 \medskip

\section{Cluster combinatorics of $d-$cluster categories}

\medskip

We now define a map $\gamma ^d_{\mathcal{H}}$ from $\mathrm{ind}
\mathcal{C}_d(\mathcal{H})$ to $\Phi ^d_{\ge-1} $. Note that any
indecomposable object $X$ of degree $i$ in
$\mathcal{C}_d(\mathcal{H})$ has the form $M[i]$ with $M\in
\mathrm{ind}\mathcal{H}$, and if $i=d$ then $M=P_j$, an
indecomposable projective representation.

\begin{defn}

Let  $\gamma^d_{\mathcal{H}}$ be defined as follows. Let $M[i]\in
\mathrm{ind}\mathcal{C}_d(\mathcal{H})$, where $M\in \mathrm{ind}H$
and $i\in \{ 1, \cdots, d\}$ (note that if $i=d$ then $M=P_j$ for
some $j$). We set
$$\gamma^d_{\mathcal{H}}(M[i])=\left\{
\begin{array}{lrl}(\underline{\mathrm{dim}}M)^{i+1} & \mbox{ if } & M[i]\in \mbox{ind}\mathcal{H}[i],
 \mbox{ for some } 0\leq i\leq d-1 ;\\
&&\\
(-\alpha _j)^1& \mbox{ if } &M[i]=P_j[d],\end{array}\right.$$

\end{defn}

This map is one kind of extension of correspondence in Gabriel-Kac's
Theorem between the indecomposable representations of quivers and
positive roots of corresponding Lie Kac-Moody algebras. It is a
bijection if $\G$ is a Dynkin diagram.

We denote by $\Phi^{sr}_{>0}$ the set of real Schur roots of $(\G,
\Omega)$, i.e.
$$\Phi^{sr}_{>0}=\{ \underline{\mathrm{dim}} M\ \ M\in \mathrm{ind} \mathcal{E(H)}\ \}.$$
  Then the map $M\mapsto \underline{\mathrm{dim}} M$ gives a 1-1
  correspondence between $\mathcal{E(H)}$ and $\Phi^{sr}_{>0}$
  \cite{Rin}.

 If we denote the set of colored almost positive real Schur roots by $\Phi^{sr, d}_{\geq -1}$, which consists of,
  by definition, $d$ copies of the set $\Phi^{sr}_{>0}$ together with one copy of the negative simple roots, then the map   $\gamma^d_{\mathcal{H}}$ gives a bijection from
  $\mathcal{E}(\mathcal{C}_d(\mathcal{H}))$ to  $\Phi^{sr, d}_{\geq -1}$.  $\Phi^{sr, d}_{\geq -1}$ contains
  a subset $\Phi^{sr,d}_{>0}$ consisting
 of all colored positive real Schur roots. The restriction of $\gamma^d_{\mathcal{H}}$
 gives a bijection from $\mathcal{E}_+(\mathcal{C}_d(\mathcal{H}))$ to  $\Phi^{sr, d}_{>0}$.

Since $ \mathcal{E}(\mathcal{H})\longrightarrow  \Phi^{sr}_{>0}: \
M\mapsto \underline{\mathrm{dim}}M$ is a bijection,
 we use $M_{\beta}$ to denote the
unique indecomposable exceptional representation in $\mathcal{H}$
whose dimension vector is $\beta$. It follows from Proposition 4.2
that $\gamma^d_{\mathcal{H}}(M_{\beta}[i]))=\beta ^{i+1}$ for any
$0\leq i\leq d-1$. We sometimes use $M_{\beta^{i+1}}$ to denote the
unique preimage of a colored almost positive real Schur root
$\beta^{i+1}$ under $\gamma^d_{\mathcal{H}}$.
\medskip

We now prepare to define a simplicial complex $\Delta ^{d,
\mathcal{H}}(\Phi)$ associated with any root system $\Phi$, which
turns out to be isomorphic to the cluster complex $\Delta
^d(\mathcal{H})$ of the $d-$cluster category
$\mathcal{C}_d(\mathcal{H})$. When $\G$ is a Dynkin graph, taking an
alternating orientation $\Omega_0$ of $\G$, this complex $\Delta
^{d, \mathcal{H}_0}(\Phi)$ is the generalized cluster complex
$\Delta ^d(\Phi)$ defined by Fomin and Reading \cite{FR2}.
\medskip

First of all, we define the "$d-$compatibility degree" on any pair
of colored almost positive real Schur roots.

 \begin{defn} For any pair of colored almost positive real Schur roots $\alpha, \beta$, the
  $d-$compatibility degree of $\alpha , \beta$ is defined as follows:
$$(\alpha ||\beta)_{d, \mathcal{H}}=\mathrm{dim}_{\mathrm{End}M_{\alpha}}(\mathrm{Ext}^1(M_{\alpha}, \oplus_{i=0}^{i=d-1}M_{\beta}[i])),$$
where
$\mathrm{dim}_{\mathrm{End}M_{\alpha}}(\mathrm{Ext}^1(M_{\alpha},
\oplus_{i=0}^{i=d-1}M_{\beta}[i]))$ denotes the length of
$\mathrm{Ext}^1(M_{\alpha}, \oplus_{i=0}^{i=d-1}M_{\beta}[i])$ as a
right $\mathrm{End}M_{\alpha}-$module. When $d>1$,
$\mathrm{End}M_{\alpha}$ is a division algebra by Proposition
\ref{basic} (3), and this length equals the dimension of
$\mathrm{Ext}^1(M_{\alpha}, \oplus_{i=0}^{i=d-1}M_{\beta}[i])$ over
the division algebra $\mathrm{End}M_{\alpha}.$

\end{defn}

\begin{rem} When $\G$ is a Dynkin diagram with trivial valuation and $\Omega_0$ is an alternating orientation of $\G$,
this compatibility degree is defined in \cite{Th}.
 When $d=1$ and $\G$ is a Dynkin diagram, we recover the classical compatibility degree defined in \cite{BMRRT, Z2}.  \end{rem}

\begin{thm}  \begin{enumerate}
\item For any pair of colored almost positive real Schur roots $\alpha , \  \beta $, \begin{enumerate} \item
$(\alpha ||\beta)_{d,\mathcal{H}}=(\sigma_{k,d}(\alpha) ||\sigma_{k,d}(\beta))_{d, s_k\mathcal{H}},$
if $k$ is a sink (or a source);
\item $(\alpha ||\beta)_{d, \mathcal{H}}=(R_{d,\Omega}(\alpha) ||R_{d,\Omega}(\beta))_{d,\mathcal{H}},$
\item $(\alpha ||\beta)_{d, \mathcal{H}}=0$ if and only if $(\beta|| \alpha)_{d,
\mathcal{H}}=0$.\end{enumerate}
\item For any almost positive real Schur root $\beta$, $((-\alpha_i)^1 ||(\beta)^l)_{d,\mathcal{H}}=0$
 if and only if  $\mathrm{  max}\{n_i(\beta), 0\}=0$ where $n_i(\beta)$ is the coefficient of $\alpha_i$ in
 the expansion of $\beta$ in terms of the simple roots $\alpha_1, \cdots, \alpha_n$.\end{enumerate}
\end{thm}
\begin{proof} \begin{enumerate}

\item Let $\alpha , \  \beta $ be two colored almost positive real Schur
roots.
\begin{enumerate}\item
We prove it for the case $k$ is a sink, the proof for source is
similar. It is easy to check
 that the following diagram is commutative:

\[ \begin{CD}
\mathrm{ind}\mathcal{C}_d(\mathcal{H}) @>\tilde{S_k^+}
>> \mathrm{ind}\mathcal{C}_d(s_k\mathcal{H})
\\
@V\gamma^d _{\mathcal{H}} VV  @VV\gamma^d  _{s_k\mathcal{H}}V  \\
\Phi^d_{\geq -1} @>\sigma_{k,d}>> \Phi^d_{\geq -1}
\end{CD} \]

Hence we have that $$\begin{array}{ll}(\sigma_{k, d}(\alpha)
||\sigma_{k,d}(\beta))_{d, s_k\mathcal{H}}&=
\mathrm{dim}_{\mbox{End}\tilde{S^+_k}(M_{\alpha})}\mbox{Ext}^1(\tilde{S^+_k}(M_{\alpha}),
 \oplus_{i=0}^{i=d-1}\tilde{S^+_k}(M_{\beta})[i])\\
&= \mathrm{dim}_{\mbox{End}{M_{\alpha}}}\mbox{Ext}^1(M_{\alpha}, \oplus_{i=0}^{i=d-1}M_{\beta}[i])\\
&=(\alpha ||\beta)_{d, \mathcal{H}}\end{array}$$

\item
As we mentioned before, the shift functor $[1]$ of
$\mathcal{C}_d(\mathcal{H})$ is an auto-equivalence.
 We now check the following diagram commutes:

\[ \begin{CD}
\mathrm{ind}\mathcal{C}_d(\mathcal{H}) @>[1]
>> \mathrm{ind} \mathcal{C}_d(\mathcal{H})
\\
@V\gamma^d _{\mathcal{H}} VV  @VV\gamma^d  _{\mathcal{H}}V  \\
\Phi^d_{\geq -1} @>R_{d,\Omega}>> \Phi^d_{\geq -1}
\end{CD} \]

By Proposition 2.2, any indecomposable object in
$\mathcal{C}_d(\mathcal{H})$ is of the form $X[i]$ with $X$ an
indecomposable representation in $\mathcal{H}$ and with $0\le i\le
d-1$, or of the form $P_j[d]$. Denote  by
$\underline{\mathrm{dim}}X=\alpha$. If $i\le d-2$, then
$R_{d,\Omega}\gamma
^d_{\mathcal{H}}(X[i])=R_{d,\Omega}((\alpha)^{i+1})=(\alpha)^{i+2}=\gamma^d
_{\mathcal{H}}[1](X[i])$. We will prove the equality for other
indecomposable objects in $\mathcal{C}_d(\mathcal{H})$. Firstly we
have that $R_{d,\Omega}\gamma
^d_{\mathcal{H}}(P_j[d-1])=R_{d,\Omega}((\underline{\mathrm{dim}}P_j)^{d})=(-\alpha_j)^1$
and $\gamma^d _{\mathcal{H}}[1](P_j[d-1])=(-\alpha _j)^1.$ Hence
$R_{d,\Omega}\gamma ^d_{\mathcal{H}}(P_j[d-1])= \gamma^d
_{\mathcal{H}}[1](P_j[d-1]).$  Secondly, for any $X[d-1]$ with $X$
not being projective, we have $\tau X \in \mathrm{ind}\mathcal{H}$.
We have that  $R_{d,\Omega}\gamma
^d_{\mathcal{H}}(X[d-1])=R_{d,\Omega}((\alpha )^{d})=(R_{\Omega}(
\alpha))^1$ and $\gamma^d _{\mathcal{H}}[1](X[d-1])=\gamma^d
_{\mathcal{H}}(X[d])= \gamma^d _{\mathcal{H}}(\tau^{-1}[d]\tau X)=
\gamma^d _{\mathcal{H}}(\tau X)=(R_{\Omega}(\alpha))^1$. The last
equality holds since $\underline{\mathrm{dim}}\tau X=
R_{\Omega}(\underline{\mathrm{dim}}X)$ (compare Section 3.1).  This
proves that $R_{d,\Omega}\gamma ^d_{\mathcal{H}}(X[d-1])=\gamma^d
_{\mathcal{H}}[1](X[d-1]).$ Similar proof for $P_j[d]$. We finish
the proof of the commutativity of the diagram.

It follows that
$$\begin{array}{ll}(R_{d,\Omega}(\alpha)
||R_{d,\Omega}(\beta))_{d,\mathcal{H}}&=
\mathrm{dim}_{\mbox{End}M_{\alpha}[1]}\mbox{Ext}^1(M_{\alpha}[1], \oplus_{i=0}^{i=d-1}M_{\beta}[1][i])\\
&= \mathrm{dim}_{\mbox{End}M_{\alpha}}\mbox{Ext}^1(M_{\alpha}, \oplus_{i=0}^{i=d-1}M_{\beta}[i])\\
&=(\alpha ||\beta)_{d,\mathcal{H}},\end{array}$$  where the second
equality follows from the fact that $[1]$ is an equivalence.
\item Let $X, Y\in \mathcal{C}_d(\mathcal{H})$ with $\mbox{Ext}^i(X,Y)=0$ for any $1\leq i\leq d.$
Then by the Calabi-Yau property of $\mathcal{C}_d(\mathcal{H})$, we
have that for any $1\leq j\leq d \  \mbox{Ext}^j(Y,X)\cong
\mbox{Ext}^{d-j+1}(X, Y)=0$. This proves $(c).$
\end{enumerate}
\item We first prove the necessity: Let $\beta$ be an almost positive real Schur root with
 $((-\alpha_i)^1 ||(\beta)^l)_{d,\mathcal{H}}=0$. If $\beta$ is a negative simple root
and $l=1$, we have easily that max$\{n_i(\beta), 0\}=0$. Now we
assume that $\beta$ is a positive real Schur root. From the
condition $((-\alpha_i)^1 ||(\beta)^l)_{d,\mathcal{H}}=0$,  we have
$\mbox{Ext}^j(P_i[d], M_{\beta^l})=0,$ i.e. $\mbox{Ext}^j(P_i[d],
M_{\beta}[l-1])=0$ , for any $1\leq j\leq d$. Since $1 \leq  l\leq
 d$, we have $1\leq j\leq d$, where $j=d+1-l$. Now we have that
$0=\mbox{Ext}^j(P_i[d], M_{\beta}[l-1])\cong \mbox{Hom}(P_i[d],
M_{\beta}[l+j-1]) \cong \mbox{Hom}(P_i,M_{\beta})$. Hence
$n_i(\beta)=\mathrm{dim}_{\mathrm{End}P_i}\mbox{Hom}(P_i,
M_{\beta})=0.$

Now we prove the other direction. Suppose that $\beta$ is an almost
positive real Schur root with $\mathrm{  max}\{0, \ n_i(\beta)\}=0.$
Firstly, if $\beta$ is the negative of a simple root, say $(-\alpha
_j)^1$, then
$$\begin{array}{ll}((-\alpha_i)^1 ||(-\alpha _j)^1)_{d,\mathcal{H}}&=
\mathrm{ dim}_{\mbox{End}(P_i[d])}\mbox{Ext}^1(P_i[d],
 \oplus_{k=0}^{k=d-1} P_j[d][k])\\
&=\mathrm{ dim}_{\mbox{End}(P_i[d])}\mbox{Ext}^1(P_i,
\oplus_{k=0}^{k=d-1} P_j[k])
\\
&= \mathrm{ dim}_{\mbox{End}(P_i[d])}\mbox{Ext}^1(P_i,
P_j[d-1]),\end{array}$$ the last equality following from Proposition
\ref{basic} (4). But $\mbox{Ext}^1(P_i, P_j[d-1])\cong \mbox{Hom}
(P_i, P_j[-1][d+1])\cong D\mbox{Hom}(P_j[-1], P_i) \cong
D\mbox{Ext}^1(P_j,P_i)=0$. This proves that $((-\alpha_i)^1
||(-\alpha _j)^1)_{d,\mathcal{H}}=0$.
 Now we assume that $\beta$ is a positive real Schur root and $l$ is a positive integer not exceeding $d$.
 We will prove that
 $((-\alpha_i)^1 ||(\beta )^l)_{d,\mathcal{H}}=0$ under the condition that $n_i(\beta)=0$.
 We can assume that $d>1$ since for $d=1$, the
 corresponding result is proved in \cite{Z1}. It follows from the condition $n_i(\beta)=0$
 that $\mbox{Hom}_{\mathcal{H}}(P_i,M_{\beta})=0$, and then  $\mbox{Hom}(P_i,M_{\beta})=0$. Hence
 $\mbox{Ext}^1(P_i[d],M_{\beta}[d-1])\cong
\mbox{Hom}(P_i,M_{\beta})=0.$  We will prove that
$\mbox{Ext}^j(P_i[d], M_{\beta^l})=0$ for $1\leq j\leq d$. Now given
such $j$, $\mbox{Ext}^j(P_i[d], M_{\beta^l})= \mbox{Ext}^j(P_i[d],
M_{\beta}[l-1])= \mbox{Ext}^1(P_i[d], M_{\beta}[l+j-2])\cong
\mbox{Ext}^1(P_i, M_{\beta}[l+j-d-2])$. Since $1\leq l\leq d, \
1\leq j\leq d$, we have $-d\le l+j-d-2\leq d-2$. Then we have that
$\mbox{Ext}^j(P_i[d], M_{\beta^l})=0$ which follows from Proposition
\ref{basic} (4) for $l+j-d-2\not=-1$ and from the fact that
$\mbox{Ext}^1(P_i, M_{\beta}[-1])\cong \mbox{Hom}(P_i, M_{\beta})=0$
for $l+j-d-2=-1$.
\end{enumerate}\end{proof}

\begin{defn} Let $\Phi$ be the root system corresponding to $\G$ and $\mathcal{H}$ the
category of representations of the valued quiver $(\G, \Omega).$
\begin{enumerate}\item Any pair  $\alpha,\ \beta$ of almost positive real Schur roots is called
 $d-$compatible if $(\alpha ||\beta)_{d, \mathcal{H}}=0$;
 a subset of $\Phi^{sr,d}_{\geq-1}$ is called $d-$compatible if any two elements of this subset
 are compatible.
\item The simplicial complex $\Delta ^{d,\mathcal{H}}(\Phi)$ associated to $\Phi$
and $\mathcal{H}$ is a complex which has $\Phi^{sr,d}_{\geq-1}$ as
the set of vertices. Its simplices are $d-$compatible subsets of
$\Phi^{sr,d}_{\geq-1}$. The subcomplex of $\Delta
^{d,\mathcal{H}}(\Phi)$ which has $\Phi ^{sr,d}_{>0}$ as the set of
vertices is denoted by $\Delta _+^{d,\mathcal{H}}(\Phi).$ We call
$\Delta ^{d,\mathcal{H}}(\Phi)$  the generalized cluster complex
associated to $\Phi$ and $\mathcal{H}$.
\end{enumerate}
\end{defn}

\begin{rem} Given a graph $\G$, we have the corresponding root system
$ \Phi$. Since the set of real Schur roots of $\Phi$ depends on the
category $\mathrm{ind}\mathcal{H}$, equivalently, on the orientation
$\Omega$ of $\G$, the generalized cluster complexes $\Delta
^{d,\mathcal{H}}(\Phi)$ are possibly non-isomorphic for different
orientations of $\G$, but they are isomorphic to each other if $\G$
is a Dynkin diagram by Proposition 4.12 (2) and the following
theorem.
\end{rem}

\begin{thm} \begin{enumerate}\item Let $\G$ be a valued graph and $\Phi$ the corresponding root system.
Let $\Omega$ be an admissible orientation of $\G$.  Then $\gamma
^d_{\mathcal{H}}$ provides an isomorphism from the simplicial
complex $\Delta ^d(\mathcal{H})$ to the generalized cluster complex
$\Delta^{d,\mathcal{H}}(\Phi),$ which sends vertices to vertices,
$k-$faces to $k-$faces.
\item The restriction of $\gamma ^d_{\mathcal{H}}$ to $\Delta ^d_+(\mathcal{H})$ gives an
isomorphism from  $\Delta ^d_+(\mathcal{H})$ to $\Delta _+^{d,\mathcal{H}}(\Phi).$
\item If $\G$ is a Dynkin graph and $\Omega _0$ is an alternating orientation of $\G$,
then $\Delta ^{d, \mathcal{H}_0}(\Phi)$ is the generalized cluster complex $\Delta^d(\Phi)$
defined by Fomin-Reading in \cite{FR2}.
\end{enumerate}
\end{thm}

\begin{proof}
\begin{enumerate} \item $\gamma^d_{\mathcal{H}}$ provides a bijection from the vertices
of $\Delta ^d(\mathcal{H})$ to that
 of $\Delta^{d,\mathcal{H}}(\Phi).$  For any pair of colored almost positive real
 Schur roots $\alpha ^k, \ \beta^l$, they are $d-$compatible if and only
 if $M_{\alpha^k}\oplus  M_{\beta^l}$ is an exceptional object where
 $M_{\alpha^k}$ and  $M_{\beta^l}$ are the  exceptional objects corresponding to
  $\alpha ^k, \ \beta^l$ respectively under the map $\gamma^d_{\mathcal{H}}$.
  Hence $\gamma^d_{\mathcal{H}}$ is an isomorphism from $\Delta ^d(\mathcal{H})$ to
  $\Delta^{d,\mathcal{H}}(\Phi).$
\item This is a direct consequence of $1.$
\item This is a direct consequence of Theorem 3.9 and Theorem 5.4.
\end{enumerate}

\end{proof}
From Theorem 5.7, one can translate results from each side. For
example, one gets the number of $d-$cluster tilting objects in
$\mathcal{C}_d(\mathcal{H})$ from the number of facets of
generalized cluster complexes  of finite root systems \cite{FR2}.

\begin{cor} \begin{enumerate}\item The generalized cluster complex
$\Delta ^{d,\mathcal{H}}(\Phi)$ and its subcomplex $\Delta
_+^{d,\mathcal{H}}(\Phi)$ are pure of dimension $n-1$.
\item Let $(\G, \Omega)$ be a connected Dynkin quiver and $\Phi$ the
 root system corresponding to $\G$. Then the number of $d-$cluster tilting
objects of $\mathcal{C}_d(\mathcal{H})$ is
$\prod_{i}\frac{dh+e_i+1}{e_i+1},$ where $h$ is the Coxeter number
of $\Phi$ and $e_1, \cdots, e_n$ the exponents of $\Phi$.
\item Let $(\G, \Omega)$ be a connected Dynkin quiver and $\Phi$ the
corresponding root system. Then the number of complements of
 any almost complete tilting object in $\mathcal{C}_d(\mathcal{H})$ is $d+1$. \end{enumerate}
\end{cor}

\begin{proof} \begin{enumerate}\item It follows from Proposition 4.12(1) and Theorem 5.7(1).
\item It follows from Theorem 5.7(1) and  Proposition 8.4 in \cite{FR2} that the statement in
2. holds for the $d-$cluster category $\mathcal{C}_d(\mathcal{H}_0)$
of $\Omega_0$. Then by Proposition 4.12(2), the statement in 2.
holds for a $d-$cluster category $\mathcal{C}_d(\mathcal{H})$
corresponding to an arbitrary orientation $\Omega$.
\item It follows from Theorem 5.7(1) and  Proposition 3.10 in \cite{FR2} that the
number of complements of
 any almost complete tilting object in $\mathcal{C}_d(\mathcal{H}_0)$ is $d+1$. Hence by Proposition
 4.12(2) the number of complements of
 any almost complete tilting object in $\mathcal{C}_d(\mathcal{H})$ is $d+1$.
\end{enumerate}\end{proof}
\begin{rem} Corollary 5.8 (1) generalizes Theorem 2.9 in \cite{FR2} to infinite root systems.\end{rem}

\begin{rem} \begin{enumerate}\item
From Corollary 5.8 (2) for $d=1$, combining with the result in
\cite{BMRRT} (see also \cite{KZ}), in which the cluster tilting
subcategories in $\mathcal{D}$ are proved
 to be in one-to-one correspondence with
the cluster tilting modules in cluster categories by the projection
$\pi,$ we have an explanation on why the number of cluster tilting
subcategories (i.e. $\mathrm{Ext}-$configurations in \cite{I1}) in
$\mathcal{D}$ is the same as the number of facets of $\Delta(\Phi)$.
\item Corollary 5.7(3) is proved by Thomas in \cite{Th} for an
alternating simply-laced Dynkin quiver $(\G,\Omega_0)$, using a
different approach .

\end{enumerate}

\end{rem}
\medskip

\begin{center}
\textbf {ACKNOWLEDGMENTS.}\end{center} The author would like to
thank Claus Michael Ringel and Osamu Iyama  for his helpful
suggestions on cluster complexes and the formalism of the paper
 and to thank Bernhard Keller for
informing the author of Y. Palu's work. After completing this work,
I was told by Bernhard Keller and
 Robert Marsh that H.Thomas had some similar results for simply-laced Dynkin type \cite{Th}. I thank them very much for this!

The author would like to thank the referees for their useful suggestions to improve the paper.%\newpage

%\newpage


\small

\begin{thebibliography}{99}


\bibitem[ABST]{ABST}
I. Assem, T. Br¨¹stle, R. Schiffler and G. Todorov,
 \newblock $m$-cluster categories and $m$-replicated algebras.
 \newblock Preprint {\tt arXiv:math.RT/0608727}, 2006.




\bibitem[AT]{AT}
C. Athanasiadis and E. Tzanaki.
\newblock Shellability and higher Cohen-Macaulay conectivity of generalized cluster complexes.
\newblock Preprint {\tt arXiv:math.CO/0606018}, 2006.


\bibitem[BaM]{BaM}
K.Baur and R.Marsh.
\newblock A Geometric description of $m-$Cluster categories.
\newblock Preprint, {\tt math.RT/0610512}. To appear in Transactions
of the AMS.

\bibitem[BM]{BM}
A.Buan and R.Marsh.
\newblock Cluster-tilting theory.
\newblock Trends in representation theory of algebras and related topics,
Edited by J.de la Pe$\tilde{n}$a and R. Bautista. Contemporary Mathematics.
 \textbf{406}, 1-30, 2006.



\bibitem[BMR]{BMR}
A.Buan, R.Marsh, and I.Reiten.
\newblock Cluster mutation via quiver representations.
\newblock Preprint, {\tt arXiv:math.RT/0412077}, 2004. To appear
in Comment. Math.Helv.


\bibitem[BMRRT]{BMRRT}
A.Buan, R.Marsh, M.Reineke, I.Reiten and G.Todorov.
\newblock Tilting theory and cluster combinatorics.
\newblock Advances in Math. \textbf{204}, 572-618, 2006.




\bibitem[CC]{CC}
P.Caldero and F.Chapoton.
\newblock Cluster algebras as Hall algebras of quiver representations.
\newblock Comment. Math. Helv. \textbf{81}, 595-616, 2006


\bibitem[CCS]{CCS}
P.Caldero, F.Chapoton and R.Schiffler.
\newblock Quivers with relations arising from clusters ($A_n$ case).
\newblock Transactions
of the AMS. \textbf{358}, 1347-1364, 2006.

\bibitem[CK]{CK}
P.Caldero and B.Keller.
\newblock From triangulated categories to cluster algebras. To appear in Invent.Math. {\tt arXiv:math.RT/0506018}.
\newblock


\bibitem[DR]{DR}
V.~Dlab and C.~M.~Ringel.
\newblock  Indecomposable representations of graphs
 and algebras.
\newblock Mem. Amer. Math. Soc. \textbf{591} (1976).

\bibitem[FR1]{FR1}
S.Fomin and N.Reading.
\newblock Root system and generalized associahedra,
\newblock lecture notes for the IAS/Park City Graduate Summer School in Geometric Combinatorics, 2004.

\bibitem[FR2]{FR2}
S.Fomin and N.Reading.
\newblock Generalized cluster complexes and Coxeter combinatorics.
\newblock IMRN. \textbf{44}, 2709-2757, 2005.

\bibitem[FZ1]{FZ1}
S.Fomin and A.Zelevinsky.
\newblock Cluster Algebras I: Foundations.
\newblock Jour. Amer. Math. Soc. \textbf{15}, no. 2, 497--529, 2002.


\bibitem[FZ2]{FZ2}
S.Fomin and A.Zelevinsky.
\newblock Y-system and generalized associahedra,
\newblock Ann of Math. \textbf{158}, 977-1018, 2003.



\bibitem[FZ3]{FZ3}
S.Fomin and A.Zelevinsky.
\newblock Cluster algebras II: Finite type classification.
\newblock Invent. Math. \textbf{154}, no.1, 63-121, 2003.



\bibitem[I1]{I1}
O.Iyama.
\newblock Higher dimensional Auslander-Reiten theory on maximal orthogonal subcategories.
\newblock To appear in Adv. Math. See also  {\tt arXiv:math.RT/0407052}.



\bibitem[I2]{I2}
O.Iyama.
\newblock Maximal orthogonal subcategories of triangulated categories satisfying Serre duality,
\newblock Mathematisches Forschungsinstitut Oberwolfach Report, no. 6(2005), 353-355.

\bibitem[IY]{IY}
O.Iyama and Y. Yoshino
\newblock Mutations in triangulated categories and rigid Cohen-Macaulay modules,
\newblock Preprint. {\tt arXiv:math.RT/0607736}.

\bibitem[Ke]{Ke}
B.Keller.
\newblock Triangulated orbit categories.
\newblock Documenta Math. \textbf{10}, 551-581, 2005.

\bibitem[KR]{KR}
B.Keller and I.Reiten.
\newblock Cluster-tilted algebras are Gorenstein and stably Calabi-Yau.
\newblock To appear in Adv.
Math. See also {\tt arXiv:math.RT/0512471}.


\bibitem[KZ]{KZ}
S.Koenig, and B.Zhu.
\newblock From triangulated categories to abelian categories-- cluster tilting in a  general framework.
\newblock To appear in Math. Zeit. See also {\tt arXiv:math.RT/0605100}.



\bibitem[MRZ]{MRZ}
R.~Marsh, M.~Reineke and A.~Zelevinsky.
\newblock Generalized associahedra via quiver representations.
\newblock Transactions of AMS. \textbf{355(10)}(2003), PP. 4171-4186.


\bibitem[Pa]{Pa}
Y. Palu.
\newblock  Ph.D Thesis in preparation.



\bibitem[Rin]{Rin}
C.~M.~Ringel. Some remarks concerning tilting modules and tilted
algebras. Origin. Relevance. Future. An appendix to the Handbook of
tilting theory, edited by Lidia Angeleri-H\"ugel, Dieter Happel and
Henning Krause. Cambridge University Press (2007), LMS Lecture Notes
Series 332.

\bibitem[Th]{Th}
H.Thomas.
\newblock Defining an $m-$cluster category.
\newblock, Preprint, 2005.

\bibitem[Z1]{Z1}
 B.Zhu.
\newblock BGP-reflection functots and cluster combinatorics,
\newblock Jour. of Pure and Applied Algebra \textbf{209}, 497-506, 2007.

\bibitem[Z2]{Z2}
B.Zhu.
\newblock Equivalences between cluster categories,
\newblock Jour. of Algebra. \textbf{304}(2006), 832-850.


\end{thebibliography}
\end{document}